\documentclass[12pt]{article}
\usepackage{amsmath}
\usepackage{amsfonts}
\usepackage{amsthm}
\usepackage{amscd}
\theoremstyle{plain}
\newtheorem{theorem}{Theorem}[section]
\newtheorem{proposition}{Proposition}[section]

\theoremstyle{definition}

\theoremstyle{remark}

\catcode `\@=11
\@addtoreset{equation}{section}

\def\section{\@startsection {section}{1}{\z@}{-3.5ex plus -1ex minus
     -.2ex}{2.3ex plus .2ex}{\normalsize\bf}}
\def\subsection{\@startsection{subsection}{2}{\z@}{-3.25ex plus -1ex minus
 -.2ex}{1.5ex plus .2ex}{\normalsize\bf}}

\def\thebibliography#1{\section*{References\markboth
  {REFERENCES}{REFERENCES}}\list
  {[\arabic{enumi}]}{\settowidth\labelwidth{[#1]}\leftmargin\labelwidth
  \advance\leftmargin\labelsep
  \usecounter{enumi}}
  \def\newblock{\hskip .11em plus .33em minus -.07em}
  \sloppy
  \sfcode`\.=1000\relax}
 
\catcode `\@=12
\def\R{\mathbb R }
\def\N{\mathbb N }
\def\P{\mathbb P }

\def\Z{\mathbb Z }
\def\C{\mathbb C }
\def\W{\mathcal{W}}
\def\V{\mathcal{V}}
\def\k{\kappa}
\def\ka{\kappa}

\def\a{\alpha}

\def\l{\lambda}
\def\w{\omega}
\def\d{\delta}

\def\del{\Delta_{\lambda}}
\def\i{{\,\mathrm i\,}}
\def\e{{\,\mathrm e}}
\def\VS{\mathrm{Vect}(S^1)}
\def\DS{\mathrm{Diff}^+(S^1)}
\def\End{\mathrm{End}}
\def\H{\mathrm{H}}
\def\deg{\mathrm{deg}}
\def\res{\mathrm{res}}
\def\ord{\mathrm{ord}}
\def\fpz{\frac {d }{dz}}
\def\pfz#1{\frac {d#1}{dz}}
\def\Fln #1{\mathcal{F}_{#1}^\lambda}
\def\cintt{\frac 1{2\pi\i}\int_{C_{\tau}}}
\def\cint{\frac 1{2\pi\i}\int_{C}}
\newcommand{\vac}[1][0]{{\vert#1 \rangle}}
\newcommand{\lpol}[1][z]{\C\,[\,#1,#1^{-1}]}
\newcommand{\RS}{Riemann surface}
\renewcommand{\L}{\mathcal{L}}
\newcommand{\Lh}{\widehat{\mathcal{L}}}
\newcommand{\g}{\mathfrak{g}}
\newcommand{\G}{\mathcal{G}}
\newcommand{\K}{\mathcal{K}}
\newcommand{\A}{\mathcal{A}}
\newcommand{\MA}{M\setminus A}
\newcommand{\gh}{\widehat{\mathfrak{g}}}
\newcommand{\Gh}{\widehat{\mathcal{G}}}
\newcommand{\Ah}{\widehat{\mathcal{A}}}
\newcommand{\Fl}[1][\lambda]{\mathcal{F}^{#1}}
\newcommand{\dzl}{dz^\lambda}
\def \nord #1{:\mkern-5mu{#1}\mkern-5mu:}  
\def\im{\text{Im\kern1.0pt }}
\def\re{\text{Re\kern1.0pt }}

\def\FA{Funktional Anal. Appl.}
\def\PL{Phys. Lett. B}
\newcommand{\refE}[1]{(\ref{E:#1})}
\newcommand{\refS}[1]{Section~\ref{S:#1}}
\newcommand{\refSS}[1]{Section~\ref{SS:#1}}
\newcommand{\refT}[1]{Theorem~\ref{T:#1}}

\begin{document}
\hspace*{\fill} Mannheimer Manuskripte 231

\hspace*{\fill} math.QA/9806032

\hspace*{\fill} May 1998
\vspace*{0.8cm}
\noindent

{ \bf SUGAWARA CONSTRUCTION FOR
\vskip 0.3cm
 HIGHER GENUS RIEMANN SURFACES
\footnote{Research partially supported by the 
Volkswagen-Stiftung
(RiP-program at Oberwolfach) and a DFG-RFBR cooperation
(436 RUS 113/276)}
\vspace{1.3cm}}

\noindent
\hspace*{1cm}
\begin{minipage}{13cm}
Martin Schlichenmaier \vspace{0.3cm}\\
Department of Mathematics and Computer Science\\
University of Mannheim\\
D-68131 Mannheim, Germany\\
E-mail:schlichenmaier@math.uni-mannheim.de
\end{minipage}

\vspace*{0.5cm}

\begin{abstract}
\noindent
By the classical genus zero Sugawara construction one obtains
from admissible  representations of affine Lie algebras
(Kac-Moody algebras of affine type) representations of the 
Virasoro algebra.
In this lecture first the classical construction is 
recalled.
Then, after giving a review on the global multi-point algebras 
of Krichever-Novikov type for compact Riemann surfaces of
arbitrary genus, the higher genus Sugawara construction is
introduced.
Finally, the lecture reports on results obtained 
in joint work with O.K. Sheinman.
We were able to show  that also
in the higher genus, multi-point situation 
one obtains from representations of the global
algebras of affine type representations of a centrally
extended algebra of meromorphic vector fields on
Riemann surfaces.
The latter algebra is the generalization of the Virasoro 
algebra to higher genus.

\vspace*{0.2cm}
{\sl  Invited lecture  at the XVI${}^{th}$  workshop
on geometric methods in physics,  Bialowieza, Poland,
June 30 -- July 6, 1997
}
\end{abstract}
\section{\hspace{-4mm}.\hspace{2mm}
INTRODUCTION}\label{S:intro}

Two-dimensional conformal field theory has played an important role
in theoretical physics for about 20 years now.
One of its origins lies in the realm of statistical
mechanics, where it is reflected by the scaling
invariance of phase transitions.
Another of its origins lies in string theory and 
two-dimensional quantum field theory.
Let me mention in this latter context the pioneering work of
A.A. Belavin, A. M. Polyakov and A.B. Zamolodchikov \cite{BPZ}.

The Sugawara construction \cite{sugawara} 
is one of the basic constructions in two-dimensional
conformal field theory.
For example,  starting from 
Wess-Zumino-Novikov-Witten (WZNW) models associated to a 
simple Lie group $G$ the modes of the energy-momentum tensor
define a representation 
of the Virasoro algebra, the centrally extended Lie algebra
of vector fields on the circle $S^1$.

This has its mathematical interpretation in the representation
theory of affine Lie algebras (Kac-Moody algebras of
untwisted affine type).
If one interprets the affine Lie algebra as Lie algebra of the
loop group $LG$, i.e. of the group of maps from $S^1$ to the group $G$,
and the algebra of vector fields on $S^1$ as Lie algebra of the group $\DS$
of orientation preserving diffeomorphisms on $S^1$
one obtains as an integrated version of the Sugawara construction 
the fact that every positive energy representation of the loop group 
 admits a projective intertwiner action of $\DS$.
For more information, general background, and references
to the original (mathematical) literature let me refer to the book of
A. Pressley and G. Segal \cite[\S 13]{PS}.
For references to the works of physicists see the review \cite{HKOC}.

The Sugawara construction above is the classical genus zero
construction in the sense that these vector fields on
$S^1$ which are finite linear combinations of Fourier modes can be 
extended to meromorphic vector fields on the projective
line $\P^1(\C)$, resp. on the sphere $S^2$, resp. on the Riemann surface
of genus zero, which  
are holomorphic outside $\ \{0,\infty\}$.

In the approach to conformal field theory on general compact Riemann
surfaces (or equivalently on smooth projective curves over $\C$)
due to A.~Tsuchiya, K.~Ueno and Y.~Yamada \cite{TUY} the Sugawara construction
is employed locally at every ``insertion point'' to
define the Knizhnik-Zamolodchikov connection for the bundle
of conformal blocks over the moduli space of curves.
With the help of this connection the authors proved factorization rules
and the Verlinde formula by studying the behaviour at the boundary 
components of the moduli space representing stable, but non-smooth curves.
See also the approach of physicists, e.g. T. Eguchi, H. Ooguri
\cite{EgOo} and D. Friedan, S. Shenker \cite{FS}.
For the connection of the Sugawara construction to geometric
quantization see for example J.-L. Brylinski, D. McLaughlin \cite{BrML}.

After the above indications about the importance of the Sugawara
construction, let me now outline the plan of the lecture.
Because the plenary lectures were supposed to be 
pedagogical in nature 
and because the main part of the audience might not be very 
familiar with conformal field theory I will first recall 
in \refS{classical} the 
definition of the Virasoro algebra,
the description how one obtains central extensions of Lie algebras, and how
the Sugawara construction works in the classical situation,
first for the special case of the Heisenberg algebra (the oscillator algebra)
and then in the case of general affine Lie algebras.
Only these  objects from conformal field theory we will need here.

The results presented are completely well-known and classical 
nowadays. Hence, I will omit nearly all references to the original
literature and instead give as reference for further reading and 
more details the books of V. Kac \cite{Kac} and  V. Kac,
 A.K. Raina~\cite{KaRa}.

Starting from the genus zero case in the above interpretation,
it is quite natural to try to
extend the situation to higher genus Riemann surfaces
using globally defined meromorphic objects.
The global objects, algebras, etc. were introduced by Krichever and
Novikov \cite{KNFA} in 1987 for the two-point situation and generalized by
me to the multi-point situation \cite{lmp}, \cite{diss}.
This (non-trivial) generalization is crucial for  applications in the theory
of  conformal blocks. 
In \refS{knalg} I recall the definition of these objects and the relevant
results.

Finally, \refS{sugawara} reports on the extension of  the 
Sugawara construction to higher genus Riemann surfaces with the help
of the Krichever-Novikov objects.
This is joint work with Oleg K. Sheinman from Moscow.
We obtain 
from representations of the global multi-point 
algebras of affine type representations of a centrally
extended algebra of meromorphic vector fields on
Riemann surfaces.
The latter algebra is the generalization of the Virasoro 
algebra to higher genus.
Details can be found in \cite{suga}, \cite{habil}.
See also this section for reference to related work, e.g. \cite{BRRW}.
At the moment we are employing the  higher genus Sugawara construction
on a global operator approach to the theory of conformal blocks
\cite{SScb}.

\section{\hspace{-4mm}.\hspace{2mm}
THE CLASSICAL GENUS ZERO SITUATION}\label{S:classical}

In this section I  will  recall the definition of the Virasoro
algebra and how the classical (genus zero) Sugawara construction works.

\subsection{\hspace{-5mm}.\hspace{2mm}The Witt Algebra}\label{SS:witt}

Let $S^1$ be the unit circle.
Denote by $\DS$ the group of orientation 
preserving diffeomorphisms $S^1\to S^1$. The space of $C^\infty$-vector fields
$\VS$ can be considered as its Lie algebra.
The exponential $\exp:\VS\to \DS$
is defined by assigning to a vector field the flow it generates.
Let me point out that the exponential map is neither locally 1:1 nor 
locally surjective \cite[p.28]{PS}. This is quite in contrast to
the finite-dimensional situation.
Inside the Lie algebra $\VS$ we can consider certain important subalgebras,
e.g. the subalgebra of analytic vector fields, or 
 the subalgebra of the vector fields which are sums of
finitely many Fourier modes.
With the latter subalgebra we will deal in the following.
It is convenient to pass to  complex-valued vector fields.
Let $\varphi$ be the coordinate on $S^1$ then the vector
fields
\begin{equation}\label{E:basf}
l_m:=(-\i)\e^{\i m\varphi}\frac {d}{d\varphi},\qquad m\in\Z
\end{equation}
correspond to the Fourier modes.
One calculates 
\begin{equation}\label{E:vir}
[\,l_n,l_m]=(m-n)\;l_{n+m}\ .
\end{equation}
The Lie algebra $\W$ generated by (finite) linear combinations of the
$l_m$ is called the {\em Witt algebra}, or the {\em Virasoro algebra
without central term}.
In the complex coordinate $z=\exp(\i\varphi)$
the elements  \refE{basf} can be considered as the restrictions of the
meromorphic vector fields
\begin{equation}\label{E:basc}
l_m:=z^{m+1}\frac {d}{dz},\qquad m\in\Z,
\end{equation}
which are holomorphic in
$\C\setminus \{0\}$.
If we take $w=1/z$ the local coordinate at $\infty\in \P^1(\C)$ we see
with $\ \frac d{dz}=-w^2\frac d{dw}\ $ that $l_m$ can be represented there
as 
\begin{equation}\label{E:basci}
l_m=-w^{1-m}\frac {d}{dw},\qquad m\in\Z\ .
\end{equation}
We can describe $\W$ as the Lie algebra of meromorphic vector fields
on $\P^1(\C)$ which have only poles at $\{0,\infty\}$.
This description will be our point of
generalization to higher genus.

Note that from \refE{basc} and \refE{basci} it follows that
the subalgebra of global holomorphic vector fields is the subalgebra
\begin{equation}
{\langle \,l_{-1},l_0,l_1\,\rangle}_{\C} \cong sl(2,\C)\ .
\end{equation}

\subsection{\hspace{-5mm}.\hspace{2mm}Central Extensions and 
the Virasoro  Algebra}

In the quantization  process one is naturally led from representations to
projective representations. In our context, let $\Phi:\W\to \End(V)$
be a realization of the vector fields as linear operators on a 
vector space $V$, e.g. the quantum Hilbert space 
of the theory. Typically  we obtain in important 
applications for the commutator in
$\End(V)$ 
\begin{equation}\label{E:repcent}
[\Phi(l_n),\Phi(l_m)]=\Phi([l_n,l_m])+\gamma(l_n,l_m)\cdot id_V\ ,
\end{equation}
with some bilinear form $\gamma:\W\times\W\to\C$.
If $\gamma=0$ then the realization will be a linear representation
of the Lie algebra $\W$. If not,
it still will be very often a
projective representation. 
This says that 
$\gamma$ fulfils additional
cocycle conditions, see \refE{lcocyc}.
Such projective representations can be described more conveniently 
as linear representations of a
central extension of $\W$. This algebra will be the
Virasoro algebra.

\medskip
In view of the next section let
 me recall the definition of a {\em central extension}
of a Lie algebra.
Let $L$ be an arbitrary Lie algebra. A bilinear map 
$\alpha: L\times L\to\C$ is called a {\em (two-)cocycle} if
\begin{equation}
\label{E:lcocyc}
\begin{gathered}
\alpha(a,b)=-\alpha(b,a)\ ,\\
\alpha([a,b],c)+\alpha([b,c],a)+\alpha([c,a],b)=0\ ,
\end{gathered}
\end{equation}
for all $a,b,c\in L$. Given such a cocycle
a (one-dimensional) central extension $\widehat{L}_\a$ 
can be defined by taking as underlying vector space
\begin{equation}
\widehat{L}_\a:=\C\oplus L,
\end{equation}
and as Lie product
\begin{equation}\label{E:liecent}
[(r,a),(s,b)]=(\alpha(a,b),[a,b]),\qquad r,s\in\C,\quad a,b\in L\ .
\end{equation}
By the cocycle condition \refE{lcocyc} $\widehat{L}_\a$ 
 is a Lie algebra. If we use the notation $t:=(1,0)$
and $\widehat{a}:=(0,a)$ we can rewrite \refE{liecent} as
\begin{equation}\label{E:centext}
[\widehat{a},\widehat{b}\,]=\widehat{[a,b]}+\alpha(a,b)\cdot t,
\quad a,b\in L\ ,\qquad
[\,t,\widehat{L}_\a]=0\ .
\end{equation}
In particular, $\C\cdot t$ is central in $\widehat{L}_\a$.

If we compare \refE{repcent} with \refE{centext} we see that 
by setting
$\Phi(\widehat{a}):=\Phi({a})$ and 
$\Phi(t):=id_V$ we obtain a linear representation of the centrally extended
algebra $\widehat{L}_\gamma$ (i.e.  for which the cocycle $\gamma$
 is used to define the extension).

Central extensions are classified up to equivalence
by  the elements of the 
Lie algebra cohomology group $\H^2(L,\C)$,
i.e. by the two-cocycles modulo coboundary.
(See \cite{HiSt} for the definition and more background information.)
Note, that by multiplying the cocycles with $r\in\C^*$ we obtain
non-equivalent but isomorphic central extensions.
For the Witt algebra one calculates
$\H^2(\W,\C)\cong \C$ (e.g. see \cite{GeFu}, \cite{KaRa}).
Hence, there are only two essentially different central extensions. 
The first one is  the trivial, or splitting  one,
which is given by $\alpha=0$. The second one
is obtained for any  $\alpha\ne 0$.
It is called the universal central extension.
With a suitable normalization it can be given as
(using the notation 
$L_n:=(0,l_n)$, $t:=(1,0)$)
\begin{equation}\label{E:virasoro}
[\,L_n,L_m\,]=(m-n)L_{n+m}+\frac {1}{12}(n^3-n)\,\delta^{m}_{-n}\,t, 
\qquad [\,t,L_n]=0,\quad\text{for all}\  n,m\in\Z\ .
\end{equation}
(As usual, $\delta^{m}_{n}$ is the Kronecker delta.)
The Lie algebra  $\V$ defined in this way 
is called the {\em Virasoro algebra}.

\subsection{\hspace{-5mm}.\hspace{2mm}The Grading}

Let me point out an important fact:
The Virasoro algebra becomes a graded algebra
by defining
\begin{equation}
\deg(L_n):=n,\qquad \deg(t):=0\ .
\end{equation}
Using the grading $\V$ decomposes into subalgebras 
\begin{gather}
\V=\V_+\oplus \V_{(0)}\oplus \V_-, \\
\V_+:=\langle L_n\mid n\in\N\,\rangle ,\quad
\V_{(0)}:=\langle L_0,t\rangle,\quad
\V_-:=\langle L_{-n}\mid n\in\N\,\rangle \ .
\end{gather}
In physics one is mainly interested in {\em highest weight representations}
of the Virasoro algebra, i.e. in representations
which are generated from a vacuum vector $\vac$ by $\V_-$,
the vacuum is annihilated by    $\V_+$, we have
$L_0\vac=h\vac$ with $h\in\C$, and $t\psi=c\cdot\psi$ for every vector
$\psi$ with  a fixed $c\in\C$.
 The number $c$ is called the {\em central charge},
the number $h$ the {\em weight} of the representation.
The operator $(-L_0)$ should be considered as energy operator.
Such kind of representations are ``positive energy'' representations,
because $(-L_0)$ is diagonalisable and its eigenvalues are bounded from
below. Note that one calculates for $\psi=L_{-n}\vac$
(for  $n>0$) 
\begin{equation}
(-L_0)\psi=-L_0L_{-n}\vac=(L_{-n}(-L_0)+[-L_0,L_{-n}])\vac=
(n-h)\psi\ .
\end{equation}

\subsection{\hspace{-5mm}.\hspace{2mm}The Classical Sugawara Construction
for the Heisenberg Algebra}\label{SS:sugheis}

Let $H$ be the infinite-dimensional {\em Heisenberg algebra}.
It is defined as the Lie algebra with underlying vector space
\begin{equation}
H=\C\,\oplus {\langle a_n\mid n\in\Z\rangle}_{\C} 
\end{equation}
and Lie structure
\begin{equation}\label{E:heisen}
[a_n,a_m]=n\,\delta_{-n}^{m},
\quad n,m\in\Z\ ,\qquad
[\C,H]=0\ .
\end{equation}
It is a central extension of the one-dimensional abelian Lie algebra.
The {\em bosonic Fock space} is defined as 
\begin{equation}\label{E:fock}
F=\C\,[x_1,x_2,\dots,x_n,\dots]\ ,
\end{equation}
the space of polynomials in infinitely many variables.
We define an action of $H$ on $F$ by setting 
\begin{equation}
\begin{gathered}\label{E:fockrep}
a_n=\frac {\partial}{\partial x_n},\quad 
a_{-n}=n\,x_n\cdot \ ,\qquad n\in\N\\
a_0=id,\qquad t=id\ .
\end{gathered}
\end{equation}
The $a_n$ with $n<0$ are ``{\em creation operators}'' in the sense
that they generate all vectors (the ``states'') from
the vacuum $\vac$,
which is given here by the constant function $1$.
 The $a_n$ with $n>0$ are 
``{\em annihilation operators}'' in the sense that if 
applied to a fixed vector  suitable often they map this vector to zero.
For this representation the {\em current} is defined as the formal sum
of operators
\begin{equation}\label{E:vircur}
J(z)\ :=\ \sum_{n\in\Z}a_nz^{-n-1},
\end{equation}
and the {\em energy-momentum tensor} as 
\begin{equation}\label{E:viremom}
T(z)\ :=\ -\frac 12 \nord{J(z)J(z)} \ ,
\end{equation}
where we ignore the colons for a moment.
Applied  to an element of the Fock space 
the current \refE{vircur} will produce a formal sum
of elements of the Fock space.
By the annihilation property,
the components will vanish  
below a certain degree in the variable $z$ 
(the lower bound for the degree
 will depend on the vector on which the current operates).
Now \refE{viremom} should be understood in the sense that one
applies $J(z)$ twice  on the vector and recollects the powers
of $z$ again:
\begin{gather}\label{E:virsug}
T(z)=-\frac 12\sum_{k\in\Z}\bigg(\sum_{l\in\Z}\nord{a_{k-l}a_l}\bigg)z^{-k-2}
=\sum_{k\in\Z}S_kz^{-k-2},
\\ \label{E:virsuga}
S_k:=-\frac 12 \sum_{l\in\Z}\nord{a_{k-l}a_l}\ .
\end{gather}
Note that $S_k$ is not a member of the algebra $H$. It is only
an operator on the Fock space. Up to now it is not even clear whether it
is a well-defined operator at all.
Let $v\in F$ be a fixed polynomial, then $a_l v=0$ for $l\gg 0$, because 
there will be a maximal index of the variables in the polynomial.
Note that for $k-l\ne -l$ (i.e. $k\ne 0$) 
the operators commute: $a_{k-l}a_l=a_la_{k-l}$, see \refE{heisen}.
Hence, in the defining  sum \refE{virsug} for the $S_k$ for $k\ne 0$ 
we can interchange the operators for $l\ll 0$ and only
finitely many terms will survive if we apply it to the fixed
vector $v$.
It follows that  $S_k$ $(k\ne 0)$ is well-defined.
If we continue to ignore
the colons further on, this would 
 not be true for $S_0=\sum_l\nord{a_{-l}a_l}$,
which can easily seen  
by applying it to the element $1\in F$.
We would obtain $(-\sum_{n<0} n)1$.
By the colons we indicate that we  use normal ordering
which is defined as follows
\begin{equation}
\nord{a_na_m}\quad :=\quad
\begin{cases}
a_na_m,   &m\ge n
\\
a_ma_n,   &m< n\ .
\end{cases}
\end{equation}
This forces the annihilation operators  to the right.
The operator $S_0$ becomes well-defined, the other $S_k$  do  not change.
By direct calculation we obtain
\begin{equation}\label{E:virsugrep}
[\,S_n,S_m\,]=(m-n)S_{n+m}+\frac 1{12}(n^3-n)\,\delta^{m}_{-n}\,id\ , 
\quad n,m\in\Z \ .
\end{equation}
By comparing \refE{virasoro} and \refE{virsugrep} we see that we 
get a representation of the Virasoro algebra with central charge 1.
If we would have been able to define the action without 
passing to normal ordering there would be no central term.
So necessarily normalization, regularization, quantization, ...
forces us to consider central extensions.

\subsection{\hspace{-5mm}.\hspace{2mm}The Sugawara Construction for
the Affine Lie Algebra}\label{SS:sugaff}

The construction can be extended to representations of affine
Lie algebras (i.e. Kac-Moody algebras of affine type). Let $\mathfrak g$ be a
finite-dimensional Lie algebra
with invariant, symmetric, non-degenerate bilinear form $(..|..)$
(i.e. $([a,b]|c)=(a|[b,c])$\ ).
Such a bilinear form exists for example for reductive Lie algebras, 
i.e. Lie algebras which are direct sums of simple and abelian ones.
Take for the simple factors the Cartan-Killing form
and for the abelian factors any non-degenerate bilinear form.

The {\em affine Lie algebra} is defined as follows. Take
as vector space
\begin{equation}\label{E:viraff}
\gh:=(\frak g\otimes\lpol)\oplus \C\cdot t_1,
\end{equation}
and denote for short $x(n):=x\otimes z^n$ then a Lie algebra
structure is defined by $(n,m\in\Z)$
\begin{equation}\label{E:viraffst}
[x(n),x(m)]=[x,y](n+m)+(x|y)\,n\,\delta^n_{-m}\cdot t_1,\qquad
[t_1,\gh]=0\ .
\end{equation}
The Lie algebra without the central term is called {\em loop algebra}
or {\em current algebra}.
Again these algebras admit  gradings defined by
$\deg (x(n)):=n$ and $\deg(t_1)=0$.

Let $V$ be an {\em admissible representation} of $\gh$. This says that
every fixed vector $v\in V$ will be annihilated by $x(n)$ for every 
$x\in\g$ if $n$ is big enough,
and that the central element $t_1$ operates as a scalar $c$.
The scalar is called the {\em level} of the representation.
Let $d=\dim\mathfrak g$, and 
let $\{u_i, i=1,\ldots, d\}$ be a
 basis of $\mathfrak g$ and $\{u^i\}$ the corresponding
dual basis with respect to $(..|..)$. 
The Casimir element of $\mathfrak g$ can be given as 
\begin{equation}
C_{\g}=\sum_{i=1}^{d} u_i u^i\ .
\end{equation}
Assume it operates as scalar on the adjoint representation,
and set  $\k$ as 
\begin{equation}\label{E:kappa}
2\k=\sum_{i=1}^{d} ad_{u_i}\circ ad_{u^i}\ .
\end{equation}
The assumption is for example fulfilled for $\g$ abelian or 
$\g$  a simple Lie algebra.
In the abelian case clearly $\k=0$.
If take the normalized Cartan-Killing form (the long
roots have the square length 2) in the simple case then 
$\k$ equals the 
{\em dual Coxeter number}. For $A_n=sl(n+1,\C)$ one obtains $\k=n+1$.
 Assume further on that the level $c$ of the representation is such 
that $c+\k\ne 0$. Now we define
\begin{equation}\label{E:affsug}
S_k:=-\frac 1{2(c+\k)}\sum_{l\in\Z}\sum_{i=1}^{d}\nord{u_i(k-l)u^i(l)}\ ,
\end{equation}
with  normal ordering 
\begin{equation}\label{E:normord}
\nord{x(n)y(m)}\quad:=\quad
\begin{cases}
x(n)y(m),&m\ge n
\\
y(m)x(n),&m< n .
\end{cases}
\end{equation}
The $S_k$ are well-defined operators on $V$ and one obtains 
\begin{theorem}\label{T:sugvir}
The mapping
\begin{equation}
L_k\quad\mapsto\quad S_k,\quad k\in\Z,\qquad\text{and} \quad t\mapsto id_V
\end{equation}
defines a representation of the Virasoro algebra with central charge
\begin{equation}
\frac {c\cdot \dim \mathfrak{g}}{c+\k}\ .
\end{equation}
\end{theorem}
A proof of this now classical theorem can be found for example in 
\cite{KaRa}. See also this book for references to the original
proofs.
Note that the Heisenberg algebra is the affinization of the 
one-dimensional abelian Lie algebra and that  
the Fock space representation of  the 
Heisenberg algebra is an admissible representation.
The Sugawara construction 
in \refSS{sugheis} can be considered  
as special case
of the general construction.
But it admits a more direct and simpler proof.

\section{\hspace{-4mm}.\hspace{2mm}THE GENERALIZED KRICHEVER-NOVIKOV ALGEBRAS}
\label{S:knalg}
%
%
\subsection{\hspace{-5mm}.\hspace{2mm}The Geometric Set-Up}

As explained in \refSS{witt}
the Virasoro algebra without central term can
be considered as the algebra of meromorphic vector fields
on $\P^1(\C)$ which have no poles outside $\{0,\infty\}$.
The algebra $\lpol$ of Laurent polynomials which appear
in the definition of the affine algebras is the algebra of  
meromorphic functions which have no poles outside $\{0,\infty\}$.
Extension to higher genus means that we should  replace $\P^1(\C)$
by an arbitrary compact Riemann surface.
This was started by Krichever and Novikov \cite{KNFA} 
for the two-point situation 
in 1987 and generalized by me 1989 \cite{lmp}, \cite{diss}
to the multi-point situation.
See also R. Dick \cite{dick} and V.A. Sadov \cite{sad}
for related work.
The multi-point situation is very important in the global
 approach to the theory
of conformal blocks.
You might think this generalization will be an easy thing to do.
Of course, writing down the definitions of the algebras is indeed easy.
One first complication (and this is not the only one) is the fact 
that classically we have graded algebras and we really need them 
for the representations. It is not possible anymore to
define a honest graded structure. We will replace the grading by a weaker
concept, an almost grading, which still will do the job.
The second task is to find an invariant geometric description of the
formal algebraic objects, like the cocycles etc.

For the following let $M$ be a compact \RS, 
$A$ a finite set of points of $M$ which is divided into two disjoint
non-empty subsets $I$ and $O$: 
$\ A=I\cup O,\  \#I=K\ge 1,\  \#O=L\ge 1$, $N=K+L$.
The elements of $I$ are the ``{\em in-points}'', the elements of
$O$ are the ``{\em out-points}''.
Let $\rho $ be a meromorphic differential,
holomorphic outside $A$, with exact pole order 1
at the points of $A$, (given) positive residues at $I$,
(given) negative residues
at $O$ (of course, obeying $\sum_{P\in M}\res_P(\rho)=0$).
By requiring $\rho$ to have purely imaginary periods it is fixed.
We fix as reference point another  point $R\in M\setminus A$
and set
\begin{equation}
u(P):=\re \int_R^P\rho\ .
\end{equation}
It is  a well-defined harmonic function.
The level lines
\begin{equation}
C_\tau=\{ P\in\MA \mid \quad u(P)=\tau\  \}, \quad \tau\in\R
\end{equation}
define a fibering of $M\setminus A$.
Every level line separates the in-points from
the out-points.
For $\tau\ll 0$ the level line $C_\tau$ is a disjoint union
of (deformed) circles around the points in $I$.
For $\tau\gg 0$ it is
a disjoint union
of (deformed) circles around the points in $O$.
This has an  interpretation in string theory.
The points in $I$ correspond to free incoming strings and the
points in $O$ to free outgoing strings. 
The variable  $\tau$ might be interpreted
as proper time of the string on the world sheet.
The interaction of the strings is given by the geometry
(genus, etc.) of the world sheet.

In the classical situation we have 
$g=0$, $I=\{0\}$, $O=\{\infty\}$
and 
\begin{equation}
\rho=\frac 1z \,dz,\qquad
u(z)=\re\int_1^z\frac 1w dw=\log|z|\ .
\end{equation}
Here the level lines are honest circles.

\subsection{\hspace{-5mm}.\hspace{2mm}The Almost Graded Algebras}

Let $\L$ be the space of meromorphic vector fields on $M$ which are
holomorphic outside $A$. Under the Lie bracket of vector fields
$\L$ becomes a Lie algebra.
Let $\K$ be the canonical bundle, i.e. the bundle whose
local sections are the local holomorphic differentials. For every
$\l\in\Z$ we consider the bundle $\K^\l:=\K^{\otimes\l}$.
Its local sections are given by the holomorphic 
$\l$-form, i.e.  by local holomorphic functions
which transform as they were products of $\l$ differentials.
It is understood that for $\l=0$ we set $\K^0=\mathcal{O}$, the
trivial bundle and for $\l<0$, 
$\K^\l:={(\K^*)}^{-\l}$ where  $\K^*$ denotes the dual bundle.\footnote
{After fixing a square root of the canonical bundle (a so-called
theta characteristic) everything can  be done for
$\l\in \frac 12\Z$. 
For simplicity we consider here only the case of integer $\l$.}
I will adopt the common practice to use the same symbol for the
bundle and its (locally free, invertible) sheaf of sections.

Let $\Fl$ be 
the vector space of global meromorphic sections of $\K^\l$ which are
holomorphic on $\MA$.
Special cases are  the differentials ($\l=1$),  the functions ($\l=0$),
and  the vector fields ($\l=-1$).
We will use $\A=\Fl[0]$ and the already introduced 
 $\L=\Fl[-1]$.
In the classical situation $\A=\lpol$ and $\L=\W$.
The Lie algebra $\L$ operates on $\Fl$ by taking the Lie derivative.
In local coordinates
the Lie derivative can be described as
\begin{equation}
L_e(g)_|=(e(z)\fpz).(g(z)\dzl)=
\left( e(z)\pfz g(z)+\l\, g(z)\pfz e(z)\right)\dzl \ .
\end{equation}
Here I used the same symbol for the section of the bundle and its
local representing function.
This operation makes $\Fl$ to a Lie module over $\L$.

It is not possible to define a graded structure. But we will be satisfied
by an {\em almost grading}. 
This says that we can decompose
$\L$ as direct sums of subspaces $\L_n$
\begin{equation}
\L=\bigoplus_{n\in\Z} \L_n, \quad
\dim \L_n<\infty
\quad\text{and}\quad
[\L_n,\L_m]\ \subseteq \bigoplus_{h=n+m}^{n+m+R}\L_{h}\ ,
\ n,m\in\Z\ ,
\end{equation}
where the integer $R$ does not depend on $n$ and $m$.
The elements of $\L_n$ are called homogeneous elements of degree $n$.
$\L$ is called an almost graded  Lie algebra.
A  similar definition works for the modules.

The grading I introduce is induced by the splitting of $A$ into $I$ and
$O$. Roughly, it is given by the zero orders of the forms at the points in $I$.
More precisely, I exhibit basis elements of $\Fl$
\begin{equation}
f_{n,p}^\l\in\Fl,\quad n\in\Z,\ p=1,\ldots,K=\# I
\end{equation}
and define the homogeneous subspace of degree $n$ to be
\begin{equation}\label{E:fln}
\Fln n:=\langle\; f_{n,p}^\l\mid p=1,\ldots, K\;\rangle.
\end{equation}
I do not want to write down the descriptions for the elements in all
cases. They can be found in \cite{diss},\cite{lmp}.
To get an idea let  me consider  the following cases.
Let $K=\#I=\#O$,  and let $g=0$ or ($g\ge 2$ and $\l\ne 0,1$). Assume that
$\  I=\{P_1,P_2,\ldots,P_K\}\ $,
$O=\{Q_1,Q_2,\ldots,Q_K\}\ $
are points which for $g\ge 2$ are in generic position.
Then there exists for every $n\in\Z$ and every $p=1,\ldots,K$ up to
multiplication with a scalar a  unique element
$f_{n,p}^\l\in\Fl$ with zero-orders
\begin{equation}
\begin{aligned}
\ord_{P_i}(f_{n,p}^\l)&=(n+1-\l)-\d_{i}^{p},\qquad i=1,\ldots,K,\\
\ord_{Q_i}(f_{n,p}^\l)&=-(n+1-\l),\qquad \qquad i=1,\ldots,K-1,\\
\ord_{Q_K}(f_{n,p}^\l)&=-(n+1-\l)+(2\l-1)(g-1)\ .
\end{aligned}
\end{equation}
After choosing local coordinates at the points in $I$ the scalar can be 
fixed. 
For the other cases the orders  are only different at the points in $O$.

The elements of the basis obey the important duality relation
(after fixing of the scalar)
\begin{equation}\label{E:kndual}
\cintt f_{n,p}^\l\cdot f_{m,r}^{1-\l}=\d_{n}^{-m}\cdot
\d_{p}^{r}\ ,
\end{equation}
where $C_\tau$ is any non-singular level line.
If we introduce the formal infinite sum
\begin{equation}\label{E:delta}
\del (Q,Q'):= \sum_{n\in\Z}\sum_{p=1}^K f_{n,p}^\l(Q) f_{-n,p}^{1-\l}(Q')
\end{equation}
we can interpret it (using the duality) as the {\sl delta distribution
of weight $\l$}.
Indeed, we obtain for $f\in\Fl$ 
\begin{equation}\label{E:int}
f(Q)=\cintt f(Q')\del(Q,Q') \ .
\end{equation}
Note that from the infinite sum \refE{delta} only finitely many terms will
contribute to the calculation of the integral \refE{int}.

Let us consider the classical situation. Here
$f_n^\l(z)=z^{n-\l}\dzl$. With respect to the (quasi)-global
coordinates $z,v$
\begin{equation}
\Delta_\l(z,v)=\sum_{n\in\Z}z^{n-\l}v^{-n-(1-\l)}\dzl dv^{1-\l}\ .
\end{equation}
For $\l=0$ we obtain
\begin{equation}
\Delta_0(z,v)=\frac 1v\sum_{n\in\Z}z^{n}v^{-n} dv\ .
\end{equation}

\subsection{\hspace{-5mm}.\hspace{2mm}Central Extensions and
Affine Algebras of Higher Genus}

We also have to deal with central extensions. The strategy for defining them
is to find  geometric descriptions of the cocycles in the classical situation
which allow
for  generalizations to higher genus Riemann surfaces.
First we consider the vector field algebra. Let $e$ and $f$ be vector fields
and let us identify them with their local representing functions.
Define
\begin{equation}\label{E:cocyc}
\chi_{C,R}(e,f)=
\frac 1{24\pi\i}\int_{C}
\left(\frac 12(e'''f-ef''')-R\cdot(e'f-ef')\right) dz
\end{equation}
Here ${}'$ means differentiation with respect to the local
coordinate $z$, $R$ is a meromorphic projective connection
which is holomorphic in $\MA$
(see \refSS{procon} for its definition), and $C$ is a closed cycle in $\MA$.
Direct calculation \cite{diss} 
shows that $\ \chi_{C,R}\ $ is indeed a two-cocycle.
For its value only the homology class of
 $C$ in $\H_1(\MA,\Z)$ is of importance.
The cohomology class of \refE{cocyc} in $\H^2(\L,\C)$, hence 
the equivalence class of the central extension
defined by the cocycle is independent of 
$R$. 
For the special cycles $C=C_\tau$ (for any $\tau$) 
the value of the integral is independent of the $\tau$ chosen.
In the following we will always take $C_\tau$ as integration cycle and
drop the cycle and the connection in the notation.
By calculating residues of the integrand in \refE{cocyc} we see
that  there exists a constant $T$, such that for homogeneous
vector fields $e$ and $f$ 
\begin{equation}
\chi(e,f)\ne 0 \quad\implies\quad T\ \le\  \deg(e)+\deg(f)\ \le\  0\ .
\end{equation}
This property is called the {\em locality} of the cocycle.
For this property $C=C_\tau$ is essential.
Using the cocycle we obtain a centrally extended vector field algebra
$\Lh=\C\oplus \L$ with
\begin{equation}
[\widehat{e},\widehat{f}\,]:=
\widehat{[e,f]}+\chi(e,f)\cdot t\ ,
\end{equation}
using the above introduced notation $\widehat{e}=(0,e)$ and $t=(1,0)$.
By setting 
$\deg(\widehat{e}):=\deg(e)$ and $\deg (t):=0$ we can extend the  grading
of $\L$ to $\Lh$.  
The locality of the cocycle implies that $\Lh$ is almost graded.

Let us now consider the (associative) algebra of functions $\A$.
It is again an almost graded algebra, i.e. it admits a decomposition
$\A=\bigoplus_{n\in\Z}\A_n$, $\dim \A_n<\infty$ 
 with  a constant $S$ such that for 
all $n,m\in\Z$
\begin{equation}\label{E:gradA}
\A_n\cdot \A_m\quad \subseteq\quad \bigoplus_{h=n+m}^{n+m+S} \A_h\ .
\end{equation}
Let $\mathfrak{g}$ be a Lie algebra and $(..|..)$ an invariant,
non-degenerate, symmetric bilinear
form on $\mathfrak{g}$ as in \refSS{sugaff} (for example, $\mathfrak{g}$ 
reductive).
Again the  space $\G=\g\otimes \A$ carries a
Lie structure. We will  need its central extension
\begin{equation}
\begin{gathered}
\Gh=(\g\otimes\A)\oplus\C\cdot t_1,
\\
[{x\otimes g},{y\otimes h}]=
{[x,y]\otimes g\cdot h}-
\big((x|y)\cint g\,dh\big)\cdot t_1\ ,
\qquad [\,t_1,\Gh]=0\ .
\end{gathered}
\end{equation}
Here I identified $x\otimes g$ with $(x\otimes g,0)$.

In the following we will only consider extensions obtained
by  integrating over a level line $C=C_\tau$.
By defining
\begin{equation}
\Gh_n:=\frak g\otimes \A_n,\quad n\in\Z,\ n\ne 0,
\qquad
\Gh_0:=(\frak g\otimes \A_0)\oplus \C \,t_1\ ,
\end{equation}
we obtain 
again an almost graded Lie algebra.

Note that for $\frak g=\C$ and $(\alpha|\beta)=\alpha \beta$ we obtain a 
central extension $\Ah$ of the abelian Lie algebra $\A$.
It is a nice and easy exercise to show that in the classical situation
everything reduces to the well-known setting introduced in 
\refS{classical}. 
By the almost-grading we obtain a decomposition
\begin{equation}
\Gh=\Gh_+\oplus \G_{(0)}\oplus \Gh_-\ ,
\end{equation}
with 
\begin{equation}
\Gh_+=\bigoplus_{n>0} \Gh_n\ ,\quad
\Gh_{(0)}=\bigoplus_{n=-S,\ldots,0} \Gh_n\ , \quad
\Gh_-=\bigoplus_{n<-S} \Gh_n,\ .
\end{equation}
The constant $S$ is the constant appearing in \refE{gradA}.
The subspaces $\Gh_+$ and $\Gh_-$ are subalgebras of
$\Gh$ and can be identified with subalgebras of $\G$.
For the two-point case the affine algebras have been introduced in
\cite{KNFA}  and its representation theory has been studied in \cite{shein}.
The multi-point generalization was introduced in \cite{claus}.

\subsection{\hspace{-5mm}.\hspace{2mm}
Definition of a Projective Connection}
\label{SS:procon}

Let $\ (U_\alpha,z_\alpha)_{\alpha\in J}\ $
 be a covering of the Riemann surface
by holomorphic coordinates, with transition functions
$z_\beta=f_{\beta\alpha}(z_\alpha)=h(z_\alpha)$.
A system of local (holomorphic, meromorphic) functions
$\ R=(R_\alpha(z_\alpha))\ $ 
is called a {\em (holomorphic, meromorphic) projective
connection} if it transforms as
\begin{equation}
R_\beta(z_\beta)\cdot (h')^2=R(z_\alpha)+S(h),\qquad\text{with}\quad
S(h)=\frac {h'''}{h'}-\frac 32\left(\frac {h''}{h'}\right)^2\ ,
\end{equation}
the Schwartzian derivative.
Here ${}'$ denotes differentiation with respect to
the coordinate $z_\a$.
Note that the difference of two projective connections
is always a usual quadratic differential.

The following result can be found in \cite{Gun}:
{\em There exists  always a holomorphic projective connection}.

\section{\hspace{-4mm}.\hspace{2mm}
THE SUGAWARA CONSTRUCTION FOR HIGHER GENUS}
\label{S:sugawara}
\subsection{\hspace{-5mm}.\hspace{2mm}The Results}

The following results have been  obtained in joint work with Oleg K. Sheinman
\cite{suga}.
Let $V$ be an admissible module of the affine algebra $\Gh$ for
higher genus. This says that the central element operates as
scalar, the level of the representation,   and that 
\begin{equation}
\text{for all}\ v\in V:\quad \Gh_nv=0 \quad\text{for}\quad n\gg 0\ .
\end{equation}
Examples of such representations are the highest weight representations
studied by Sheinman \cite{shein}.

Again we have to choose a normal ordering $\nord{...}$ similar to
\refE{normord}.
For $x\in\mathfrak{g}$ we define the associated current 
$\widetilde{x}(Q)$ to be
\begin{equation}
\widetilde{x}(Q)=\sum_{n\in\Z}\sum_{p=1}^K x(n,p)\,\w^{n,p}(Q)\ .
\end{equation}
Here I used the notation
\begin{equation}
e_{n,p}:=f_{n,p}^{-1},\quad
a_{n,p}:=f_{n,p}^{0},\quad
\w^{n,p}:=f_{-n,p}^{1},\quad
\Omega^{n,p}:=f_{-n,p}^{2}\ ,
\end{equation}
(note the inverted index  for $\w^{n,p}$ and $\Omega^{n,p}$)
and $ x(n,p)$ for the operator 
corresponding to $x\otimes a_{n,p}$.
In the following, sums over the indices of the basis elements will
occur. To simplify the formulas I will drop the range of the summation
if it is clear from the situation. In particular,
the summation for the first index is always over $\Z$, for the
second from $1$ to $K$. 
Concerning the Lie algebra $\g$ 
recall the notations introduced in \refSS{sugaff}.

The {\em Sugawara field} or {\em energy-momentum field} is defined as
\begin{equation}
T(Q):=\frac 12\sum_{i=1}^{\dim\mathfrak {g}}\nord{\widetilde{u_i}(Q)
\widetilde{u^i}(Q)}\ .
\end{equation}
If we plug in the definition of the currents we obtain
\begin{equation}
T(Q)=\frac 12\sum_{n,m,p,r}\sum_i\nord{u_i(n,p)u^i(m,r)}\w^{n,p}(Q)
\w^{m,r}(Q)  \ .
\end{equation}
If $T(Q)$ is well-defined at all it has to be a form of weight 2
(a quadratic differential) in $Q$. Hence we can write it as the sum
\begin{equation}
T(Q)=\sum_{k,s}L_{k,s}\Omega^{k,s}(Q),
\end{equation}
 with 
certain operators  $L_{k,s}$, $k\in\Z$, $s=1,\ldots, K$
 on the representation space.
To calculate the operators
we use duality \refE{kndual}
\begin{equation}
L_{k,s}=\cintt T(Q)e_{k,s}(Q)=
\sum_{n,m,p,r}\sum_i\nord{u_i(n,p)u^i(m,r)}l_{(k,s)}^{(n,p),(m,r)}
\end{equation}
with
\begin{equation}
l_{(k,s)}^{(n,p),(m,r)}=\cintt 
\w^{n,p}\w^{m,r}e_{k,s}\ .
\end{equation}
By considering the poles of the  integrand we see 
for fixed $k$ that $l_{(k,s)}^{(n,p),(m,r)}\ne 0$ 
is possible only for  a finite range of values of the sum $m+n$.
E.g. for $\#I=\#O=1$ one obtains as range: $k\le m+n\le k+g$.
Altogether we still obtain infinitely many non-vanishing coefficients.
But by the normal ordering the operators with index $m\gg 0$
or $n\gg 0$ (corresponding to  $m\ll 0$) will be moved to
the right and a fixed vector $v$ will be annihilated by them.
Applied to a fixed $v$ only finitely many terms will survive in the result.
 Hence,
\begin{proposition}
The operators $L_{k,s}$ are well-defined linear operators on 
an admissible representation space $V$.
\end{proposition}
\noindent
Note that these operators are not elements of the affine algebra. 

The key result is  the following
\begin{theorem}\label{T:semi}
\cite[Prop. 3.2, Prop. A.1]{suga}
Let $c$ be the level of the representation and let $\k$ be defined by
\refE{kappa} then
\begin{equation}\label{E:ad}
[L_{k,s},x(n,p)] = -(c+\ka)x(\nabla_{e_{k,s}}a_{n,p})\ .
\end{equation}
\end{theorem}
By $x(\nabla_{e_{k,s}}a_{n,p})$ we understand the following:
take the Lie derivative of the function $a_{n,p}$ with respect to
the vector field $e_{k,s}$, tensorize the resulting  function with the element
 $x\in\mathfrak {g}$ and take the corresponding operator on $V$.
\begin{theorem}\label{T:suga}
\cite[Thm. 3.1, Thm. A.1]{suga}
Let $\mathfrak {g}$ be a finite-dimensional Lie algebra such that the 
Casimir element of $\mathfrak {g}$ acts as $2\kappa\cdot id$ on the 
adjoint representation (which is for example the case for
$\g$ abelian or simple). 
Let $V$ be an admissible representation of level $c$ of the 
affine algebra $\Gh$. Assume $c+\ka\ne 0$. Set
\begin{equation}\label{E:sugop}
L_{k,s}^*:=\frac {-1}{c+\ka}L_{k,s}
\end{equation}
then the map
\begin{equation}
\widehat{e}_{k,s}\quad\mapsto\quad L_{k,s}^*,\qquad k\in\Z,\ s=1,\ldots,K
\end{equation}
defines a representation of a  central extension of the 
(Krichever-Novikov) vector field algebra $\Lh$ given by the above 
introduced local geometric cocycle \refE{cocyc} (with a suitable projective
connection) with central charge
\begin{equation}
c_{\Lh}=\frac {c\cdot\dim\mathfrak {g}}{c+\ka}\ .
\end{equation} 
(\ $\widehat{e}_{k,s}$ is the lift of the vector field ${e}_{k,s}$).
\end{theorem}
For the proof see \cite{suga}, \cite{habil}, resp. the comments below.
The construction does depend on the normal ordering in the
sense that a different choice as \refE{normord}
would change the defining cocycle for the central extension. 
But in \cite{suga} it is  shown that 
a different normal ordering does not change the cohomology class of the
cocycle. Hence, we obtain equivalent central extensions.

\refT{suga} specializes to \refT{sugvir}
in the classical situation.
In addition, it gives a  geometric proof for it.
For $g\ge 1$ and two points the abelian (Heisenberg algebra) case 
 was considered by
Krichever and Novikov \cite{KNFA}. For the case
 that $\mathfrak {g}$ is simple
and two points the situation was studied by the physicists
Bonora, Rinaldi, Rosso and Wu \cite{BRRW}.
There the presentation was not so that every mathematician might 
accept it. So we reproved it for this situation, 
insuring that their result was correct, obtained finer
results and generalized it to the multi-point situation.
We felt us inspired by \cite{BRRW}.

\subsection{\hspace{-5mm}.\hspace{2mm}Few Comments on the Proofs}

Here I will only point out some ideas of the proofs.
For notational reason I
consider  the two point case and drop everywhere the second index.
Choose a  cut-off function $\psi$ such that 
$\psi(x) =1$ for $|x|\le 1$ and  
$\psi(x) =0$ for $|x|> 1$, $x\in\R$.
We set
\begin{equation}
L_{k}(\epsilon)=
\sum_{n,m}\sum_i\nord{u_i(n)u^i(m)}l_{k}^{n,m}\psi(\epsilon n)\ .
\end{equation}
If we apply $L_k$ to a fixed $v\in V$ 
only finitely many terms in the sum will contribute. 
Hence,  we get $L_kv=L_k(\epsilon)v$ 
for $|\epsilon|$ small enough. 
The advantage is that as long as $\epsilon\ne 0$ we can  for the
evaluation of the commutator
\begin{equation}\label{E:commu}
[L_k(\epsilon),x(r)]
\end{equation}
forget about the normal ordering.
Now
\begin{equation}
[u_i(n)u^i(m),x(r)]=u_i(n)[u^i(m),x(r)]+
[u_i(n),x(r)]u^i(m)
\end{equation}
resolves the commutators in \refE{commu}.
This gives terms connected with the level $c$ coming from the action
of the central element
and terms of the type
\begin{equation}
\sum_{n,m,s}\sum_i u_i(n)[u^i,x](s)\,\alpha_{m,r}^s l_k^{n,m}\psi(\epsilon n)\ ,
\end{equation}
where  the $\alpha_{m,r}^s$ 
are defined by $a_m\cdot a_r=\sum_s \alpha_{m,r}^s\,a_s$.
For $\epsilon\to 0$ the split infinite sums 
 considered separately do not make sense.
Before we can let $\epsilon\to 0$ 
we have to rearrange things and pass over to normal ordering again.
This gives additional terms of the type
\begin{equation}
\sum_i[u_i,[u^i,x]]  \ .
\end{equation}
And this introduces $\ka$ as half of the eigenvalue of the 
Casimir operator in the adjoint representation.
We also make essential use of the duality and the 
properties of the delta distribution $\Delta_\l$.

To prove \refT{suga} we proceed in a similar manner. Again 
consider $[L_n(\epsilon),L_m]$.
 By the cut-off we can forget about the  normal ordering,
do the calculation (using 
\refT{semi}, duality, etc.) in a non-straight-forward 
manner.
Changing back to normal ordering to make things well-defined for 
$\epsilon\to 0$ will give contributions to the central term. 

The details can be found in \cite{suga} and \cite{habil}.

\subsection{\hspace{-5mm}.\hspace{2mm}Further Results}

{\bf (a)} 
For the higher genus algebras 
by the almost-grading it is possible to introduce
a notion of highest weight representations.
In particular, the highest weight representation 
is an admissible representation.
In the classical situation the weight 
for an affine algebra is an element
of the dual of the Cartan subalgebra.
In the higher genus 
the weight has more components.
For the two-point situation this was developed by O.K. Sheinman
\cite{shein}. There you can also find the classification of the
representations.
For the multi-point case see \cite{habil}.
If we start with highest weight representations of the
affine algebra $\Gh$ we obtain highest weight representations
for the vector field algebra.
In the article \cite{suga} 
relations between the involved weights are given.

\bigskip
\noindent
{\bf (b)} 
For admissible representations  \refT{semi}
allows to consider the semi-direct product of the
operator algebra corresponding to $\Gh$ and the
algebra obtained by the Sugawara operators 
$L_{k,s}^*$ \refE{sugop}
using the rescaled action \refE{ad}
of $L_{k,s}^*$ on $x(n,p)$.
The representation of the affine algebra will naturally extend
to a representation of the 
semi-direct product.

\bigskip
\noindent
{\bf (c)} 
Under additional assumptions on the representations
it is possible to construct Casimir operators for  representations
of the affine algebra \cite{suga}.



\end{document}